# Enhanced Representative Days and System States Modeling for Energy Storage Investment Analysis

Diego A. Tejada-Arango , *Member, IEEE*, Maya Domeshek, Sonja Wogrin , *Member, IEEE*, and Efraim Centeno

*Abstract*—This paper analyzes different models for evaluating investments in energy storage systems (ESS) in power systems with high penetration of renewable energy sources. First of all, two methodologies proposed in the literature are extended to consider ESS investment: a unit commitment model that uses the "system states" (SS) method of representing time; and another one that uses a "representative periods" (RP) method. Besides, this paper proposes two new models that improve the previous ones without a significant increase of computation time. The enhanced models are the "system states reduced frequency matrix" model which addresses short-term energy storage more approximately than the SS method to reduce the number of constraints in the problem, and the "representative periods with transition matrix and cluster indices" (RP-TM&CI) model which guarantees some continuity between representative periods, e.g., days, and introduces long-term storage into a model originally designed only for the short term. All these models are compared using an hourly unit commitment model as benchmark. While both system state models provide an excellent representation of long-term storage, their representation of short-term storage is frequently unrealistic. The RP-TM&CI model, on the other hand, succeeds in approximating both short- and long-term storage, which leads to almost 10 times lower error in storage investment results in comparison to the other models analyzed.

## Nomenclature

In the following formulation "*p/s*" refer to the parameters used to identify time divisions: periods (e.g., 1 h) in the detailed model and states in the system states model respectively.

### A. Indices and Sets

| | |
|---|---|
| $p \in P$ | Periods (hours). |
| $p_l(p)$ | Subset with the last period of the time horizon. |
| $s, s' \in S$ | System states. |
| $k \in K$ | Periods in which storage limit constraints are imposed in system states models. |
| $g \in G$ | Generation units (thermal or storage). |
| $t(g)$ | Subset of thermal generation units. |
| $h(g)$ | Subset of storage units. |
| $h_l(g)$ | Subset of long-term storage (e.g., hydro) units. |
| $h_s(g)$ | Subset of short-term storage (e.g., batteries) units. |
| $n, n' \in N$ | Electrical nodes o buses. |
| $n_s(n)$ | Subset of electrical nodes or buses without slack bus. |
| $c$ | Circuits. |
| $\mathcal{G}_{gn}$ | Generators $g$ connected to bus $n$. |
| $\Theta_{nn'c}$ | Circuits $c$ connected between bus $n'$ and $n$. |
| $rp \in RP$ | Set of representative periods (e.g., days, weeks). |
| $\Gamma_{rpp}$ | Injective map of each period $p$ to a representative period $rp$. |
| $H_{pp'}$ | Injective map of each period $p$ to a period $p' \in \Gamma_{rpp}$. |
| $p_f(p, rp)$ | Subset with the first period $p$ of the representative period $rp$. |

### B. Parameters

| | |
|---|---|
| $C_g^{fuel}$ | Cost of consumed fuel [k€/MJ]. |
| $\alpha_g$ | Variable term of fuel consumption [MJ/GWh]. |
| $\beta_g$ | Fixed term of fuel consumption [MJ]. |
| $\gamma_g$ | Fuel consumption during the startup [MJ]. |
| $C_g^{om}$ | Cost of operation and maintenance [k€/GWh]. |
| $D_{p/sn}$ | Electricity demand per node [GW]. |
| $V_{p/sn}^{max}$ | Renewable production per node (e.g., wind or solar) [GW]. |
| $Q_g^{max}, Q_g^{min}$ | Upper and lower bound on production [GW]. |
| $SRR_g$ | Maximum 10-minute ramp [GW]. |
| $X^{res}$ | Operating reserve [p.u.]. |
| $W0_h$ | Initial storage level [GWh]. |
| $W_h^{max}, W_h^{min}$ | Upper and lower bound on energy storage [GWh]. |
| $W_h^{fin}$ | Minimum final storage level [GWh]. |
| $I_{p/sh}$ | Hourly energy inflows [GWh]. |
| $\eta_h$ | Efficiency of storage unit [p.u.]. |
| $B_h^{max}$ | Upper bound on charging/pumping [GW]. |
| $T_s$ | Duration of state [h]. |
| $TC_{nn'c}^{max}$ | Transmission capacity of circuit $c$ [GW]. |
| $ISF_{nn'c\,n_s}$ | Injection Shift Factors [p.u.]. |
| $N_{s\,s'}$ | Transition matrix between states. |
| $F_{ss'k}$ | Frequency matrix between states and changes. |

Manuscript received September 19, 2017; revised December 21, 2017 and February 16, 2018; accepted March 16, 2018. Date of publication April 9, 2018; date of current version October 18, 2018. This work was supported by Project Grant ENE2016-79517-R, awarded by the Spanish Ministerio de Economía y Competitividad. Paper no. TPWRS-01450-2017. *(Corresponding author: Diego A. Tejada-Arango.)*

D. A. Tejada-Arango, S. Wogrin, and E. Centeno are with the Escuela Técnica Superior de Ingeniería ICAI, Instituto de Investigación Tecnológica, Universidad Pontificia Comillas, Madrid 28015, Spain (e-mail: diego.tejada@iit.comillas.edu; sonja.wogrin@comillas.edu; Efraim.Centeno@iit.comillas.edu).

M. Domeshek is with the Smith College, Northampton, MA 01063 USA (e-mail: mdomeshek@smith.edu).

Color versions of one or more of the figures in this paper are available online at http://ieeexplore.ieee.org.

Digital Object Identifier 10.1109/TPWRS.2018.2819578

| | |
|---|---|
| $RFM_{ss'k}$ | Reduced Frequency Matrix between states and changes. |
| $WG_{rp}$ | Weight of representative periods [h]. |
| $NRP_{rp\,rp'}$ | Transition matrix between representative periods. |
| $NP_{rp}$ | Number of periods at each representative period [h]. |
| $M$ | Moving window for storage level [h]. |
| $C_h^{inv}$ | Investment cost for storage units [k€/GW]. |
| $EPR_h^{\max}$ $EPR_h^{\min}$ | Maximum and minimum energy to power ratio [h]. |

### C. Variables

| | |
|---|---|
| $q_{p/sg}$ | Power production [GW]. |
| $\hat{q}_{p/sg}$ | Power production above $Q_g^{\min}$ [GW]. |
| $v_{p/sn}$ | Renewable production [GW]. |
| $r_{p/sg}$ | Spinning reserve [GW]. |
| $w_{p/sh}$ | Storage level [GWh]. |
| $\Delta w_{ss'h}$ | Difference in storage [GWh]. |
| $b_{p/sh}$ | Hourly charged/pumped power [GW]. |
| $sp_{p/sh}$ | Hourly energy spillage [GWh]. |
| $pf_{p/s\,nn'c}$ | Power flow per circuit [GW]. |
| $pns_{p/sn}$ | Power not supply per node [GW]. |
| $u_{p/sg}$ | Binary dispatch decision [0-1]. |
| $y_{p/sg}$ | Binary startup decision [0-1]. |
| $y_{s\,s'g}$ | Binary startup decision for state model [0-1]. |
| $x_h$ | Storage investment [GW]. |

## I. INTRODUCTION

### A. Motivation

AMONG the different power system planning models, there are short-term models with high time resolution such as unit commitment models, with information pertaining to every hour, half hour, or 10 minutes; and long-term models such as investment models that ignore small time-scale changes so as to make the calculations in a reasonable amount of time.

The introduction of variable renewable energy sources (RES) into the energy system, however, makes it necessary to include more short-term dynamics, such as varying wind or sunlight availability, in long-term models [1]. Models that incorporate information at both time scales include the TIMES modeling framework [2], the Regional Energy Deployment System (ReEDS) framework [3], and the Resource Planning Model (RPM) [4]. These models have multi-year investment decisions as well as 'time slices' within each year that represent a wide variety of possible demand and RES production levels.

The time slices structure allows the models to find solutions on a representative set of situations that the system operator must be able to respond to. However, while they do not include in detail every hour of the time horizon, the calculations are not overly burdensome.

Nowadays, energy storage systems (ESS) have become a promising flexible option to deal with the variability of renewable energy [5]. Realistically modeling ESS requires the preservation of chronological information, because the amount of stored energy available at any given moment depends on the amount of energy stored in all previous time periods [6]. Although some models have endeavored to incorporate ESS investment decisions, they do not preserve chronological information and so do not fully model storage evolution [7], [8]. In this paper We created medium and long-term optimization models for ESS investment with reduced representation of time that nevertheless maintains some chronology for the sake of co-optimizing different types of storage technologies. Moreover, we propose some new models to improve the existing ones in the literature.

### B. Literature Review

There are two common ways to reduce temporal information while maintaining some chronology that can be found in the literature: '*representative periods*' and '*system states*'. The system states are also referred to as load periods, load duration curves, or time slices in more simplified versions. Both methods are based on clustering techniques. In this section, we describe the main characteristics of both methods and review publications that present them.

In the '*representative periods*'[1] (RP) method, a certain number of days, groups of days, or in some cases weeks that are representative of the variety of situations that can be found during the course of the time horizon (e.g., year) are chosen. All calculations (e.g., investment decisions and unit dispatch) are done for the selected days or weeks. Each RP 'represents' the periods in the year that are similar to itself, so one can reconstruct the behavior of the system over the whole year by using the values calculated for the RPs in place of the periods they represent. The RPs preserve the internal chronology of their hours, making for a more realistic representation of changing storage level over the course of a day or week. However, the RP method does not preserve the chronology among the RPs. Therefore, any ESS with a cycle[2] longer than the RP (e.g., weekly monthly, or yearly rather than daily) will not be chronologically represented with the highest accuracy. This method has been used for some of the models that try to incorporate both long- and short-term dynamics, such as the RPM model in [9]. There has been much debate about the best way to choose these RPs. Some authors use a heuristic method, choosing one day for each season or one day in each season for the week and for the weekend. Others have proposed methods that involve optimizing both the number and clustering of RPs to minimize the difference between the load duration curve and the approximate one created by the RPs [10], [11]. There has also been debate about the optimal length for RPs. For instance in [12], the authors suggested representative groups of days or representative weeks, whose advantage is that it increases the amount of chronology preserved, and whose disadvantage is, of course, that it increases the calculational burden. The most versatile method for grouping RPs comes from [13], and relies on clustering techniques (e.g., k-means or k-medoids) to group a number of hours with any number of normalized characteristics (solar energy, demand, wind energy, etc). No matter how long the periods or how they are chosen, the drawback of the RP method is that it can only deal with relatively short-term storage cycles, those that charge and discharges in the course of a period (e.g., day), but not, for example, with hydro reservoirs with monthly or yearly cycles.

---

[1]In this paper, we use the name 'representative periods' when general concepts and model formulation are explained. However, for the case study and results, we use the name 'representative days' because the selected period is equal to a day.

[2]One cycle here refers to the total amount of time that it takes the ESS to go from minimum operating capacity to maximum operating capacity and back to minimum capacity.

The other method, *'system states'* (SS) was introduced in [14]. It is designed to be an improvement on the entirely non-chronological load duration curve method. The SS method characterizes each time step (e.g., hours) in the time horizon by a set of features such as demand, wind, and solar power availability. Hours with similar values of these features are considered to belong to the same *'system state'*. Every hour in the time horizon is then assigned to one of the system states, and calculations are done for each system state in the same way they would be done for each hour of an hourly model. As with the representative periods, each system state gets a weight or duration that depends on the number of real time periods in the time horizon that are represented by it. This is also called time slices in models such as ReEDS [15]. The innovation of SS method in [14] is the transition matrix, which counts up the number of transitions between all system states, allowing the addition of chronological constraints, such as start-up constraints. In [16] the system states method was extended to deal with storage. Although each system state can only calculate the change in energy storage, the total storage in any given hour can be calculated ex post by adding up all the changes in storage from the beginning of the time horizon to the hour of interest. The total storage is kept within bounds during the modeling process by backtracking to calculate the total storage at certain chosen hours in the time horizon and constraining storage in those hours to be in bounds. This idea was applied and analyzed in [17] for the operation of a network-constrained power system. This paper further extends the use of this SS method to the ESS investment problem.

As we mention at the beginning of this section, this work focuses on the reduction of temporal information. However, there are other types of reduction techniques to deal with the computational burden in long-term planning models, such as transmission network aggregation [18], [19]; exogenous estimation of curtailment reduction, curtailment itself, and capacity value [9], [20]. These methods are compatible with the models proposed in this paper and could be combined to further improve the reduction of the computational burden. Nevertheless, these sorts of combinations are beyond the scope of this work.

### C. Contributions

The first aim of this paper is to compare the SS method and the RP method for an ESS investment model in order to determine which one is better or what system characteristics the quality of the approximation method depends on. However, we found some difficulties and drawbacks in the basic formulation of both methods, which are explained in Section III.D. Therefore, the second aim of this paper is to develop enhanced versions of both methods in order deal with these difficulties. Thus, the main contributions of this paper are:

1) The extension of the SS method in [16], [17] to consider ESS investment.
2) The formulation of enhanced versions of SS and RP to preserve the chronological information of different kinds of ESS cycles (from hourly to yearly), which outperform existing methods in terms of solution quality and CPU time and allow for the co-optimization of both short- and long-term storage.
3) The comparison of SS and RP for ESS investment models using an hourly unit commitment model as a benchmark.

The paper is organized as follows: Section III shows model formulations used for SS and RP, including the proposed enhanced formulations for both methods. Section IV analyzes the results in a Spanish case study based on European visions for the year 2030. Section V discusses the benefits of considering a unified modeling approach with different operating profiles (e.g., seasonal and intraday). Finally, Section V concludes this paper.

## II. MODEL FORMULATION

This section contains the five model formulations compared in this paper.

### A. Hourly Unit Commitment Model (HM)

The following equations describe the hourly unit commitment model used as the benchmark to test the proposed models, which is based on [21].

$$\min_\Omega \sum_{p,t} \left\{ C_t^{fuel} \cdot [\beta_t u_{pt} + \gamma_t y_{pt} + \alpha_t q_{pt}] + C_t^{om} q_{pt} \right\}$$
$$+ \sum_h C_h^{inv} x_h. \quad (1a)$$

Subject to:

$$\sum_{t \in \mathcal{G}} q_{pt} + \sum_{h \in \mathcal{G}} (q_{ph} - b_{ph}) + v_{pn} + \sum_{n'c \in \Theta} (pf_{pn'nc} - pf_{pnn'c})$$
$$+ pns_{pn} = D_{pn} \quad \forall\, p, n \quad (1b)$$

$$pf_{pnn'c} = \sum_{n_s} ISF_{nn'cn_s} \cdot \left[ \sum_{t \in \mathcal{G}_{tn_s}} q_{pt} + \sum_{h \in \mathcal{G}_{hn_s}} (q_{ph} - b_{ph}) \right.$$
$$\left. + v_{pn_s} + pns_{pn_s} - D_{pn_s} \right]$$
$$\forall\, nn'c \in \Theta, p \quad (1c)$$

$$q_{pt} = Q_t^{\min} u_{pt} + \hat{q}_{pt} \quad \forall\, p, t \quad (1d)$$

$$0 \leq \hat{q}_{pt} \leq (Q_t^{\max} - Q_t^{\min}) u_{pt} \quad \forall\, p, t \quad (1e)$$

$$u_{pt} - u_{p-1,t} \leq y_{pt} \quad \forall\, p, t \quad (1f)$$

$$r_{pt} + q_{pt} \leq u_{pt} Q_t^{\max} \quad \forall\, p, t \quad (1g)$$

$$0 \leq r_{pt} \leq SRR_t \quad \forall\, p, t \quad (1h)$$

$$\sum_t r_{pt} \geq X^{res} \cdot \sum_n D_{pn} \quad \forall\, p \quad (1i)$$

$$u_{pt}, y_{pt} \in \{0,1\} \quad \forall\, p, t \quad (1j)$$

$$w_{ph} = w_{p-1,h} + W0_{p=1,h} + I_{ph} - q_{ph} - sp_{ph} + \eta_h b_{ph} \quad \forall\, p, h \tag{1k}$$

$$0 \leq v_{pn} \leq V_{pn}^{\max} \quad \forall\, p, n \tag{1l}$$

$$0 \leq q_{ph} \leq Q_h^{\max} + x_h \quad \forall\, p, h \tag{1m}$$

$$0 \leq b_{ph} \leq B_h^{\max} + \eta_h x_h \quad \forall\, p, h \tag{1n}$$

$$0 \leq sp_{ph} \quad \forall\, p, h \tag{1o}$$

$$|pf_{pnn'c}| \leq TC_{nn'c}^{\max} \quad \forall\, nn'c \in \Theta, p \tag{1p}$$

$$W_h^{\min} + EPR_h^{\min} x_h \leq w_{ph} \leq W_h^{\max} + EPR_h^{\max} x_h \quad \forall\, p, h \tag{1q}$$

$$w_{p_l, h} \geq W_h^{fin} \quad \forall\, h \tag{1r}$$

The objective function (1a) minimizes storage investment costs and the total operating cost of the system (e.g., startup costs, fixed costs, variable costs, operations and maintenance costs, and penalties for spillage and energy not supplied). Constraint (1b) is the demand balance equation. Constraint (1c) represents the power flow equation using Injection Shift Factors (ISF). Constraints (1d)–(1e) ensure thermal unit production is within minimum and maximum capacity. (1f) is the startup constraint of the unit-commitment. Equations (1g)–(1i) are reserve constraints. Equation (1j) states that the commitment and connection variables are binary. Equation (1k) is the storage constraint which states that the storage in any hour is the storage in the previous hour plus the net charging and discharging in the current hour. Equations (1l)–(1q) keep within bounds the renewable production per node, the power output per storage unit, the pumped power per storage unit, the energy spillage, the power flow through a line, and the amount of energy stored in each storage unit. Equations (1m) and (1n) include the power capacity increase due to the storage investment variable. Equation (1p) includes the energy capacity increase considering parameters $EPR_h^{\max}$ and $EPR_h^{\min}$. These parameters describe the relationship between the energy that can be stored (maximum and minimum respectively) and the nominal power of the equipment. Finally, constraint (1r) establishes the minimum storage level at the last period of the time horizon.

### B. System States Model (SS)

This section presents the formulation of the system states model as conceived in [17].

$$\min_{\Omega} \sum_{s,t} \left\{ C_t^{fuel} \cdot \left[ T_s \beta_t u_{st} + \sum_{s' \neq s} N_{s's} \gamma_t y_{s'st} + T_s \alpha_t q_{st} \right] \right. \\ \left. + C_t^{om} T_s q_{st} \right\} + \sum_h C_h^{inv} x_h \tag{2a}$$

$$\sum_{t \in \mathcal{G}} q_{st} + \sum_{h \in \mathcal{G}} (q_{sh} - b_{sh}) + v_{sn} + \sum_{n'c \in \Theta} (pf_{sn'nc} - pf_{snn'c}) \\ + pns_{sn} = D_{sn} \quad \forall\, s, n \tag{2b}$$

$$pf_{snn'c} = \sum_{n_s} ISF_{nn'cn_s} \cdot \left[ \sum_{t \in \mathcal{G}_{tn_s}} q_{st} + \sum_{h \in \mathcal{G}_{hn_s}} (q_{sh} - b_{sh}) \right. \\ \left. + v_{sn_s} + pns_{sn_s} - D_{sn_s} \right] \\ \forall\, nn'c \in \Theta, s \tag{2c}$$

$$q_{st} = Q_t^{\min} u_{st} + \hat{q}_{st} \quad \forall\, s, t \tag{2d}$$

$$0 \leq \hat{q}_{st} \leq (Q_t^{\max} - Q_t^{\min}) u_{st} \quad \forall\, s, t \tag{2e}$$

$$u_{st} - u_{s',t} \leq y_{s'st} \quad \forall\, s, t \tag{2f}$$

$$r_{st} + q_{st} \leq u_{st} Q_t^{\max} \quad \forall\, s, t \tag{2g}$$

$$0 \leq r_{st} \leq SRR_t \quad \forall\, s, t \tag{2h}$$

$$\sum_t r_{st} \geq X^{res} \cdot \sum_n D_{sn} \quad \forall\, s \tag{2i}$$

$$u_{st}, y_{s'st} \in \{0, 1\} \quad \forall\, s, t \tag{2j}$$

$$0 \leq v_{sn} \leq V_{sn}^{\max} \quad \forall\, s, n \tag{2k}$$

$$0 \leq q_{sh} \leq Q_h^{\max} + x_h \quad \forall\, s, h \tag{2l}$$

$$0 \leq b_{sh} \leq B_h^{\max} + \eta_h x_h \quad \forall\, s, h \tag{2m}$$

$$0 \leq sp_{sh} \quad \forall\, s, h \tag{2n}$$

$$|pf_{snn'c}| \leq TC_{nn'c}^{\max} \quad \forall\, nn'c \in \Theta, s \tag{2o}$$

$$\Delta w_{ss'h} = 0.5 \cdot (I_{sh} + I_{s'h} + \eta_h b_{sh} + \eta_h b_{s'h} - q_{sh} - q_{s'h} \\ - sp_{sh} - sp_{s'h}) \quad \forall\, s, s', h \tag{2p}$$

$$\sum_{\substack{s,s'\ s.t.\\ N_{ss'} > 0}} N_{ss'} \cdot \Delta w_{ss'h} \geq W_h^{fin} - W0_h + EPR_h^{\min} x_h \quad \forall\, h \tag{2q}$$

$$\sum_{\substack{s,s'\ s.t.\\ N_{ss'} > 0}} N_{ss'} \cdot \Delta w_{ss'h} \leq W_h^{\max} - W0_h + EPR_h^{\max} x_h \quad \forall\, h \tag{2r}$$

$$\sum_{\substack{s,s'\ s.t.\\ F_{ss'k} > 0}} F_{ss'k} \cdot \Delta w_{ss'h} \geq W_h^{\min} - W0_h + EPR_h^{\min} x_h \quad \forall\, h, k \tag{2s}$$

$$\sum_{\substack{s,s'\ s.t.\\ F_{ss'k} > 0}} F_{ss'k} \cdot \Delta w_{ss'h} \leq W_h^{\max} - W0_h + EPR_h^{\max} x_h \quad \forall\, h, k \tag{2t}$$

The objective function (2a) incorporates storage investment and operational costs just as in the hourly model. The costs of each state are weighted by the number of hours in the time horizon that belong to that state, and the startup costs are multiplied by the transition matrix which gives the number of transitions between each set of states. Constraints (2b) to (2o) are formulated exactly as in the hourly model in Section III.A, except that they are defined for each system state 's' rather than each hour 'p'. (2p-2t) are the system states formulation of the storage constraints. Equation (2p) defines the variable $\Delta w$ which is the central difference of the net energy storage gained in two states between which there is a transition. Equations (2q) and (2r) ensure that storage in the first and last hours of the time horizon are within upper and lower bounds including the storage investment. The amount of storage in the last hour of the time horizon is determined by multiplying each $\Delta w$ by the corresponding value in the transition matrix and adding them all up. Equations (2s) and (2t) try to keep the energy storage within bounds throughout the time horizon including the storage investment. At each of the hours, $k$, a subset of all hours in the time horizon, (2s) an (2t) add up all $\Delta w$ from the beginning of the time horizon with the aid of the frequency matrices and make certain they are between maximum and minimum storage values.

### C. Representative Periods (Days/Weeks) Model (RP)

This section describes the RP model which is a commonly used method of reducing temporal information. Although the model is general enough to work with RPs of any length, we will speak of representative days for the sake of simplicity. The formulation is roughly the same as that of the hourly model, except the constraints only apply to the hours within the representative days.

$$\min_\Omega \sum_{p,rp\epsilon\Gamma_{rpp}} \left\{ WG_{rp} \cdot \sum_t \left\{ C_t^{fuel} \cdot [\beta_t u_{pt} + \gamma_t y_{pt} + \alpha_t q_{pt}] + C_t^{om} q_{pt} \right\} \right\} + \sum_h C_h^{inv} x_h \quad (3a)$$

Subject to:
Equations (1b) – (1r) $\forall\, p\epsilon\Gamma_{rpp}$

$$w_{p=p_f(p,rp)+NP_{rp}-1,h} \geq w_{p=p_f(p,rp),h} \quad \forall\, (p,rp)\,\epsilon\Gamma_{rpp}, h \quad (3b)$$

The objective function (3a) minimizes the storage investment cost and operational cost just as in the hourly model, except that the operational costs associated with each day are multiplied by the number of days in the time horizon that are represented by it to yield the cost for the entire time horizon. The RP model is constrained to (1b) to (1r) from the HM benchmark model. Nevertheless, in the RP model, (1b) to (1r) only apply to hours belonging to the selected representative days.

Equation (3b) is a special constraint introduced into the RP model that guarantees that the amount of energy stored in each unit at the end of each representative day is greater than or equal to the amount of energy in storage at the beginning of the day. Since each day is calculated separately, this prevents a unit from finishing a day with less energy than the starting level of the next day, and thus creating energy from nothing. This is a very simple way to deal with the maximum energy storage per year. Other approaches ensure that the change accumulated over each representative period does not exceed the storage limits, and ensure balance over the whole year. However, for the sake of simplicity, these types of approaches are not analyzed in this paper.

Despite the incorporation of (3b), each representative day is independent of the others and the RP model does not guarantee chronological continuity among the representative days for the ESS.

### D. Comments About System States and Representative Periods Models

The SS and RP models have some drawbacks, which are detailed in a case study in Section IV. In this section, we summarize these drawbacks:
- The SS model results and CPU time are highly dependent on (2s) and (2t). These equations guarantee that storage levels are between the maximum and minimum for each storage unit throughout the time horizon and help to keep some chronological information in the optimization process. Equations (2s) and (2t) do, however, have two disadvantages. First, short-term storage devices such as batteries require several bounds in a day to ensure that storage levels are within bounds, but the greater the number of bounds, the longer the CPU time. Second, in order to determine the number of bounds (i.e., set k size) we need an iterative process detailed in [17] which adds even more CPU time to the SS model.
- The RP model solves each representative period (e.g., day) independently and with the same constraints as the HM model. CPU time thus depends on the number of representative periods instead of on the number of bounds for storage units, as it does in the SS Model. The main drawback is that chronology among the representative periods is lost and storage levels of storage units with a cycle longer than the representative period (e.g., hydro units) are not determined adequately. This is especially important in hydrothermal power systems or power systems with pumped hydro storage potential.

In the following sections, we propose enhanced versions of the SS and RP models to tackle these drawbacks.

### E. System States Model With Reduced Frequency Matrix (SS-RFM)

This section shows the formulation of the System States Reduced Frequency Matrix Model, hereafter SS-RFM. This is a new variation on the system states model created to reduce the computational time and avoid the iterative process for determining storage bounds constraints.

Objective function: Equation (2a)
Subject to:
Equations (2b)–(2r)

$$\sum_{\substack{s,s'\,s.t.\\F_{ss'k}>0}} F_{ss'k}\cdot\Delta w_{ss'h_l} \geq W_{h_l}^{\min} - W0_{h_l} + EPR_h^{\min} x_h \quad \forall\, h_l,k \quad (4a)$$

$$\sum_{\substack{s,s'\,s.t.\\F_{ss'k}>0}} F_{ss'k}\cdot\Delta w_{ss'h_l} \leq W_{h_l}^{\max} - W0_{h_l} + EPR_h^{\max} x_h \quad \forall\, h_l,k \quad (4b)$$

$$\sum_{\substack{s,s'\,s.t.\\RFM_{ss'k}>0}} RFM_{ss'k}\cdot\Delta w_{ss'h_s} \geq W_{h_s}^{\min} - W0_{h_s} + EPR_h^{\min} x_h \quad \forall\, h_s,k \quad (4c)$$

$$\sum_{\substack{s,s' \text{ s.t.} \\ RFM_{ss'k} > 0}} RFM_{ss'k} \cdot \Delta w_{ss'h_s} \leq W_{h_s}^{\max} - W0_{h_s}$$
$$+ EPR_h^{\max} x_h \quad \forall h_s, k \qquad (4b)$$

The objective function (2a) and constraints (2b)–(2r) are exactly the same as in the SS model. The difference between the two models lies in the handling of storage which has been separated into long- and short-term storage, each with its own set of constraints. Equations (4a) and (4b) take the same form as (2s) and (2t), but are only applied to long-term storage, which is likely to go through only one or two cycles per year. Set $k$ is a subset of hours in the time horizon in which the upper and lower bound are checked. At each hour $k$, (2s) and (2t) use the frequency matrices to add up all changes in storage from the beginning of the time horizon to hour $k$ and check that the total is within bounds. Equations (4c) and (4d), represent the storage constraints for short-term storage. At each hour $k$, they add up all the net changes in storage since the last hour $k$ and constrain that sum to be within bounds. This is done with the aid of the Reduced Frequency Matrix (RFM), an innovation of this model which is just the difference between the frequency matrix ($F_{ss'k}$) corresponding to the current hour $k$ and that corresponding to the previous element in set $k$, that is, $k-1$. In other words, the difference between these two elements or hours in the set $k$ could be understood as a moving window. It is important to mention that despite the use of the RFM, the storage level could be out of bounds because the hours in set $k$ are predefined in the model and we do not know in advance the storage level value at each hour in set $k$. The best practice for reducing the number of hours in which the storage levels can be out of bounds is to predefine the moving window considering the smallest storage cycle in the power system.

### F. Representative Periods Model With Transition Matrix and Cluster Indices (RP-TM&CI)

This section shows the Representative Period with Transition Matrix and Cluster Indices (RP-TM&CI) model which is the second original contribution of this paper. Although the model is sufficiently general to be able to work with representative periods of any length, we will once again speak of representative days for the sake of simplicity.

Objective function: (3a)
Subject to:
Equations (1b)–(1r) $\forall p \epsilon \Gamma_{rpp}$

$$u_{p'=p_f(p',rp')+NP_{rp}-1,t} = u_{p=p_f(p,rp),t}$$
$$\forall t, (p, rp) \epsilon \Gamma_{rpp}, rp' / NRP_{rp\ rp'} > 0 \qquad (5a)$$

$$w_{ph} = w_{p-M,h} + W0_{p=1,h} + \sum_{p'=p-M+1}^{p} \sum_{p'' \in H_{p'p''}}$$
$$\times (I_{p''h} - q_{p''h} - sp_{p''h} + \eta_h b_{p''h}) \quad \forall p, h \qquad (5b)$$

The objective function has the same formulation as the regular representative day model, i.e., (3a). The RP-TM&CI model is constrained with (1b) to (1r) for all the hours belonging to the selected representative days. Equation (5a) is an innovation of this model. It creates continuity between the representative days and prevents unnecessary startups by using a transition matrix to require that for any pair of representative days that transition from one to the other, the thermal units that are on in the last hour of the first are also on in the first hour of the second. As written here, if there is even one transition between the two days, this constraint is applied. However, the constraint could be set to take effect only if there is a considerable number of transitions between the two days, 5 or 10% of the transitions in the time horizon, for example. Equation (5b) is the second innovation of this model; it creates the continuity in storage across the entire time horizon that allows for the modeling of long-term storage. It does this by checking at regular intervals (1 week) that all the energy charged and discharged since the previous week plus the total energy at the last check point are within bounds. This is possible because, as a result of the clustering procedure to determine the representative days, we know the Cluster Indices (CI), which is a numeric column vector where each row indicates the cluster assignment (i.e., representative day) of the corresponding day of the year. This information is included in the model using the subset $H_{pp'}$.

### III. CASE STUDIES AND RESULTS

As a case study, we chose the Spanish power system in target year 2030. The Spanish case is interesting because it has hydro reservoirs (i.e., ESS with monthly or yearly cycle) and, according to ENTSO-E [22], the next ten years will likely bring investment in Battery Energy Storage System (BESS) and Pumped Hydroelectric Energy Storage (PHES), i.e., ESS with daily or weekly cycle. We ran four different scenarios or visions for 2030 on the hourly model and the four approximate models. The wind and solar profiles for these visions were taken from [23], [24] while hourly demand data and annual production per technology were taken from the ENTSO-E *'Ten Year Network Development Plan 2016'* [22]. Vision 1 and 3 were based on national predictions, whereas visions 2 and 4 were designed with the whole of Europe and climate protection goals in mind. The scenarios include a significant development of renewable electricity sources, supplying 35% to 60% of the total annual demand, depending on the Vision. Moreover, the hourly demand curve of each Vision reflects the potential for demand response, which rises from 5% in Vision 1 to 20% in Vision 4. A summary of the main assumptions of each vision can be found in the Appendix.

For each of the four visions, the SS and RP models were run with four different numbers of clusters for increasing time resolution. The RP and RP-TM&CI models used 4, 9, 18, and 37 representative days which corresponds respectively to 1%, 2%, 5% and 10% of the time horizon. Time resolution within each representative day is hourly. The SS and SS-RFM models used 26, 48, 96, and 216 system states. These numbers of states were chosen because they provided a *'fair'* comparison with the clusters used with the RP models by having roughly the same number of binary variables.

The representative days were chosen by normalizing time series for the hourly demand, wind availability, solar availability, and hydro inflows, and combining 24 hours of those time series

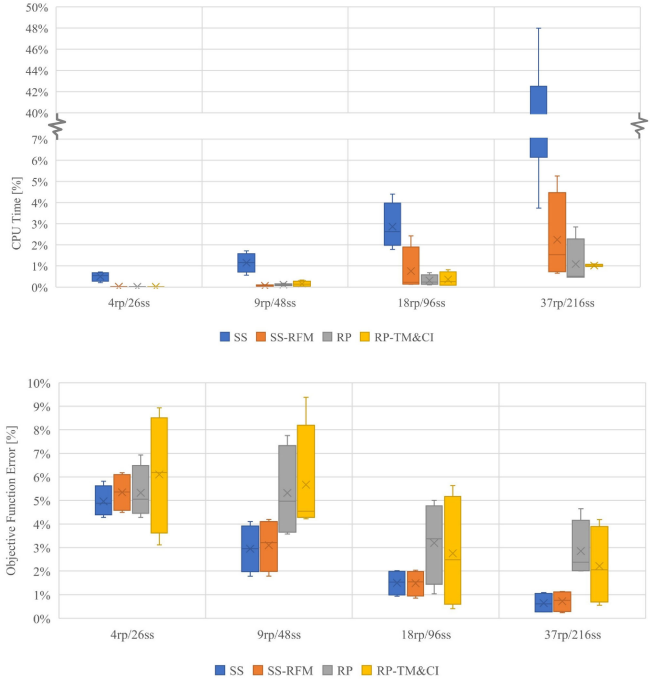

Fig. 1. CPU time (top). Objective function error (bottom).

(96 dimensions in all) into a single point to be clustered with the rest of days of the year using k-medoids. The system states were chosen in an analogous manner. The four-time series were normalized, but this time each point to be clustered represented only one hour (4 dimensions) and the clustering method was k-means so that the resulting system state was the centroid of the cluster (a composite hour) rather than a true hour.

We performed two analyses. In the first one, we ran the models without ESS investment in order to determine the accuracy of the models from the operational point of view. In the second, we analyzed the ESS investment to compare the results of investment decisions made by the four approximate models to those of the benchmark, HM model.

### A. Operation Only Results

For this case study, we considered a total BESS installed capacity of 10 GWh with a maximum output of 1 GW and a 0.9 efficiency coefficient.

Fig. 1 shows a box & whisker plot for CPU Time and objective function error considering the results for each vision. All models were solved until optimality, i.e., until the integrality gap equaled zero. Fig. 1 (top) shows the time necessary for the solution of each model as a fraction of the time taken by the hourly model as the number of clusters (i.e., system states or representative days) increases. As expected, the amount of time necessary for model solution increases with the temporal resolution, but up to the 3rd time resolution (18rp, 96ss) all four approximate models took less than 5% of the time that the hourly models took. Also, as expected, increasing the number of system states or representative days reduced the error in the objective function, see Fig. 1 (bottom). Fig. 1 also shows the improvement obtained with the SS-RFM and RP-TM&CI models proposed in this paper. The SS-RFM model took between 4 and 20 times less CPU time than the SS model without hampering the performance of the approximation in the objective function error. Moreover, the RP-TM&CI model reduced the objective function error of RP model as the number of representative days increase without a significant rise in the CPU time. These results show some of the advantages of the model proposed in this paper. For the sake of simplicity, the rest of this section shows only the results for the 3rd time resolution (18rp, 96ss) because it has a good trade-off between CPU time and objective function error.

TABLE I
AVERAGE ERRORS

|  | Result | SS | SS-RFM | RP | RP-TM&CI |
|---|---|---|---|---|---|
| Production | Nuclear | -0.3% | -0.2% | 5.4% | -0.2% |
|  | Coal | 1.9% | 1.2% | 10.5% | -2.0% |
|  | CCGT | 2.3% | 2.8% | -10.6% | 1.3% |
|  | Hydro | -0.2% | -0.2% | -10.4% | 0.8% |
|  | Battery | 7.3% | 11.3% | -17.0% | -4.8% |
|  | Renewable | -0.5% | -0.5% | -0.4% | -0.5% |
|  | RES curtailment | 24.7% | 24.9% | 18.4% | 18.6% |
| Start-up | Coal | -53.9% | -54.3% | -52.4% | -9.3% |
|  | CCGT | -73.6% | -75.2% | -91.3% | -21.0% |
| Price | Average | -0.5% | 0.03% | 8.0% | 0.7% |
|  | Max | -25.4% | -8.5% | -22.7% | 2.1% |
|  | Min | 0.0% | 0.0% | 0.0% | 0.0% |

So far, we have used objective function error to judge the accuracy of the approximate models, nevertheless, results such as annual production per technology, total number of startups, and energy prices allow for a more detailed comparison. Table I shows the average error for these results when comparing each approximate model to the hourly model. Negative values in Table I show overestimation in the approximate model while positive values are underestimation. For thermal production SS, SS-RFM, and RP-TM&CI models have errors lower than 3% while the RP model has error between 5% and 11% because it solves each representative day individually. The SS and SS-RFM models give the estimation of total hydro production closest to that of the hourly model while the RP model gives a very poor estimate. This is because the RP model constrains the storage at the end of each day to be higher than at the beginning so hydro storage cannot evolve according to its natural yearly cycle. The RP-TM&CI model, however, does succeed in estimating the annual hydro production, which is what it was designed to do. The SS and SS-RFM models do not approximate the annual battery production very well, as the models cannot keep the energy fully within bounds throughout the time horizon. The RP-TM&CI model gives a value of the total annual battery production that is closest to the HM model. RES production is estimated with good accuracy (i.e., errors less than 0.5%) for all models, while the RES curtailment has more error and is underestimated in all models. However, representative periods-type models have slightly better accuracy than system states-type models. The RP model overestimates the number of necessary

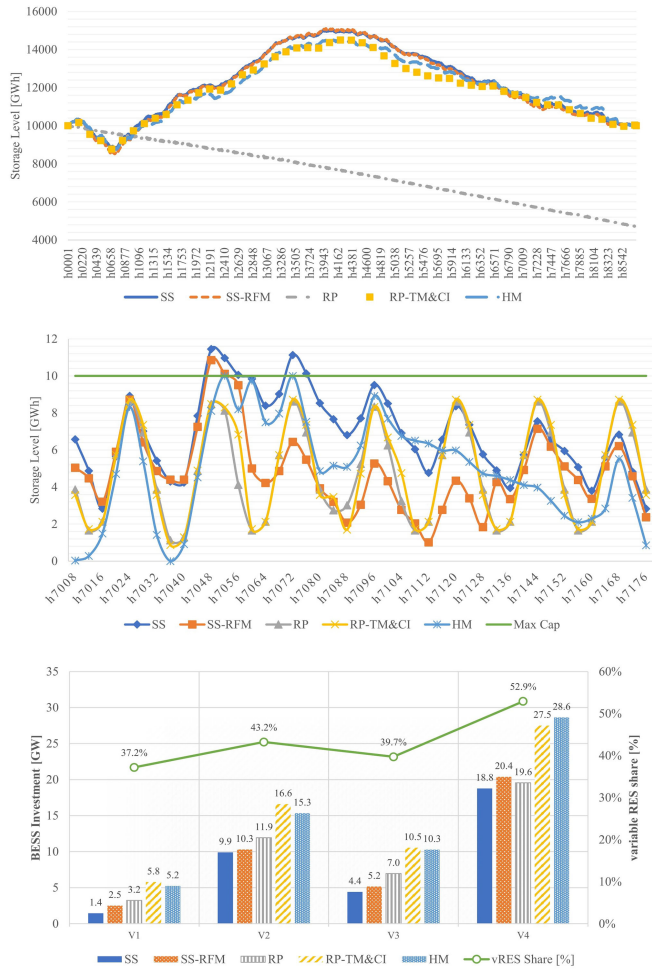

Fig. 2. Hydro storage level (top). BESS storage level (middle). BESS investment and variable RES share for each vision (bottom).

TABLE II
INVESTMENT RESULT ERROR PER VISION

| Result | Vision | SS | SS-RFM | RP | RP-TM&CI |
|---|---|---|---|---|---|
| Objective Function Error [%] | V1 | 0.5% | 0.1% | 0.8% | 0.1% |
| | V2 | 1.2% | 1.0% | 4.4% | 0.7% |
| | V3 | 0.5% | 0.4% | 4.8% | 5.4% |
| | V4 | 6.4% | 6.5% | 1.8% | 5.6% |
| Battery Investment Error [%] | V1 | 72.4% | 52.0% | 38.3% | -10.3% |
| | V2 | 35.4% | 32.9% | 22.2% | -8.3% |
| | V3 | 57.1% | 49.8% | 32.3% | -2.5% |
| | V4 | 34.4% | 28.8% | 31.7% | 3.9% |

RP model cannot correctly estimate the evolution of storage levels considering the production, consumption, inflows, and spillages for each representative day because the representative days are not related among themselves. The RP-TM&CI model fixes this by considering chronology among the representative days using the transition matrix and cluster indices. In fact, the RP-TM&CI model yields the prediction of hydro storage levels that is most similar to that of the HM model. The BESS storage level is shown in Fig. 2 (middle) for a week of the year. RP and RP-TM&CI models perform best when the BESS charge and discharge in a single day. If, however, the true BESS charges and discharges over the course of more than one day then the RP and RP-TM&CI have trouble approximating that, as they are limited to the representative days. Despite this, the RP-TM&CI model performs better than the RP model due to the chronological information shared among the representative days. The SS and SS-RFM models have better performance than the representative days models because they are not limited to the period length, i.e., 24 hours, and this allows them to capture charging and discharging periods longer than a day. However, as mentioned in Section III-E, the SSs model cannot guarantee that BESS storage levels stay within bounds. In Fig. 2 (middle) both SS and SS-RFM predict that BESS storage levels will exceed the upper bound, which is unrealistic in a power system operation. To correct that behavior, the number of constraints should be increased, but this vastly increases CPU time in the SS model and increases the error in the SS-RFM model. If the extra constrained hours are chosen using the iterative method, this increases the CPU time still further.

### B. ESS Investment Results

For this case study only the investment results are shown because the trend is similar to that of the operational results (e.g., production, number of start-ups, prices), see Section IV-A.

We consider the possibility of investment in BESS technology. Unlike the previous case study, BESS initial capacity is not predefined. We consider an investment cost of 20 [€/kW] for BESS according to the report "Technology development roadmap towards 2030" [25] and a maximum energy to power ratio ($EPR_h^{\max}$) of 4 hours. Table II shows objective function error and investment error for each vision using the HM model results as a reference. All four models underestimate the objec-startups during the year of peaking units (CCGT), which is only to be expected since it treats each day as separate from the others. Because they maintain some chronology between periods using the transition matrix, SS and SS-RFM do a better job of estimating startups than the RP. However, the RP-TM&CI model has the number of startups closest to that of the HM model, as it uses its transition matrix to keep continuity between the thermal units at the end of one day and the beginning of the next. These results also demonstrate the effectiveness of the RP-TM&CI model over the RP model. In the case of the energy prices, the RP model makes the worst estimate due to the previous results. The average prices in SS, SS-RFM, and RP-TM&CI are all quite accurate, but the maximum price is better estimated in the enhanced models, SS-RFM and RP-TM&CI. This is important because the storage investment results are partially driven by the differences between the maximum and minimum prices. We analyze this situation in Section IV-B.

Fig. 2 shows the storage level evolution for hydro unit and BESS for vision 1. Not only is the total yearly hydro production estimated by SS, SS-RFM, and RP-TM&CI very close to that of the HM as shown in Table I, but the overall storage evolution closely follows that of the HM, Fig. 2 (top). The tive function, especially when there is a high share of variable RES (vision 4). However, the range of the error values remains similar to those shown in Fig. 1. As for the investment error, the RP-TM&CI model offers the best approximation. This is because it is the model that most accurately estimates energy prices and energy production of each technology (Table I). Both the SS-RFM and RP-TM&CI models, the original contributions of this paper, represent significant improvements on their former versions SS and RP.

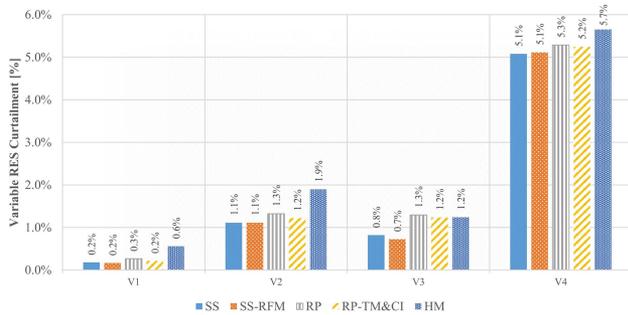

Fig. 3. Variable RES curtailment.

Fig. 2 (bottom) shows BESS investment obtained with all the models for each vision, and the share of variable RES (i.e., wind and solar productions). As expected, BESS investment increases when the variable RES share increases in the power system. The SS model and the SS-RFM model underestimate the investment by the greatest amount due to their main drawback, which is that they do not fully guarantee that the energy stored in the batteries is lower than the capacity of the batteries. This means that they permit energy to be stored beyond what investment has paid for, and therefore require less investment to achieve the same results as the RP model and the RP-TM&CI model.

## IV. Discussion

In this section we want to highlight two main aspects of the results: the relationship between RES curtailment and storage investment, and the link between short- and long-term storage.

First, Fig. 3 shows the variable RES curtailment as a percentage of the total available RES for each vision. The amount of curtailment determined by all models underestimates the reference values from the hourly model. While a portion of the under-investment in storage shown in Section IV-B is due to the inaccuracies in the way storage is represented in each model, some of the underinvestment may also come from the models' underestimation of variable RES curtailments. This is based on the tight connection between RES curtailment and storage needs, as shown in [20]. Models such as ReEDS and RPM use exogenous estimations to relate these two aspects in systems with high share of RES. However, the models proposed in this paper determine this relationship endogenously. Improvements in the clustering process could be performed to improve this relationship; however, further research is needed to verify this hypothesis.

Second, in this paper we focus on modeling energy storage investment with operational detail, considering long-term (i.e., seasonal) hydro storage generation as well as short-term (i.e., hours) storage systems such as batteries. These are very different resources in the power system. Therefore, the following question arises: Why try to model both with the same methodology? Hydro storage already exists in most real power systems and more could be built in the future, and short-term storage (e.g., BESS) is getting cheaper and could be a good technical solution to reduce RES curtailments even with relatively low energy to power ratios (e.g., 1-4 h). Moreover, if both types of storage are not considered at the same time, then an assumption must be included regarding storage operation. For example, it is possible to consider maximum available hydro energy without tracking the storage level, or to assume a peak shaving for short-term ESS. In either case one decision is fixed while the other is optimized. Therefore, possible synergies between both storage systems are neglected. This is the case of more traditional hydrothermal dispatch models.

The RP-TM&CI model co-optimizes both types of storage. Hence, the operational decisions of short- and long-term storage are now linked and depend on each other. The benefits of this co-optimization are shown in the results of Section IV. In fact, the best results are obtained with the RP-TM&CI model, which represents the relationship between both types of storage better than the other approximate models. It should also be noted that the RP-TM&CI model could be used to improve traditional hydrothermal models in which the water value serves as a consistent way of coupling long-term reservoir management with short-term operations of storage units. Using the RP-TM&CI model it might be possible to obtain the water value of long-term reservoirs internalizing the information of short-term storage, which is not possible in traditional hydrothermal models.

## V. Conclusion

This paper compares four different methods of approximating time representations in an hourly unit-commitment model with ESS investment. These methods include the SS model and the RP method as well as enhanced versions of the SS and RP models (the SS-RFM model and the RP-TM&CI models) which are the new contributions of this paper and perform better than the original versions.

The SS model was originally developed to include chronology and high time resolution details in mid- and long-term models. While it can deal with long-term storage, it cannot accurately estimate short-term storage, and quickly becomes calculation intensive because of the storage constraints. The SS-RFM model takes much less time to run than the regular SS model, because it reformulates the storage constraints, but it does not improve the accuracy of the short-term storage modeling. Moreover, SS models could lead to infeasible results (i.e., more energy stored than the maximum storage capacity), which is their major drawback, and means that they require additional adjustments for most practical applications.

Unlike the SS models, the RP model cannot handle long-term storage, but it deals well with short-term storage as it preserves within-day chronology. The RP-TM&CI model combines aspects of the SS and RP models to account for both short and long-term storage. According to the case study results, it is the most accurate of the four approximate models and does not require a significant increase of CPU time. These results support the idea that including chronological information among representative periods may be an efficient way to include small time scale variations in longer-term planning models that involve storage. Doing so is a critical need in the adequate representation of power systems that include a significant and increasing quota of variable renewable sources and energy storage systems.

TABLE III
INSTALLED CAPACITY PER VISION

| Installed capacity (MW) | Vision 1 | Vision 2 | Vision 3 | Vision 4 |
| --- | --- | --- | --- | --- |
| Gas | 24948 | 21572 | 29208 | 29208 |
| Hard coal | 5900 | 5900 | 4160 | 4160 |
| Hydro | 23450 | 23450 | 25050 | 25635 |
| Nuclear | 7120 | 7120 | 7120 | 7120 |
| Others non-RES | 10480 | 10480 | 12210 | 12210 |
| Others RES | 2400 | 2400 | 5100 | 5100 |
| Solar | 16800 | 33150 | 25000 | 54130 |
| Wind | 35750 | 27650 | 39300 | 40604 |

TABLE IV
ANNUAL GENERATION PER VISION

| Annual generation [GWh] | Vision 1 | Vision 2 | Vision 3 | Vision 4 |
| --- | --- | --- | --- | --- |
| Others non-renewable | 46438 | 46438 | 54103 | 54103 |
| Others renewable | 12587 | 12587 | 26748 | 26748 |
| Wind | 78223 | 60291 | 86414 | 89032 |
| Solar | 39313 | 69870 | 58266 | 112707 |
| Run-of-river | 19814 | 19814 | 19814 | 19814 |
| Hydro generation | 15697 | 16651 | 16118 | 20392 |
| Nuclear | 49943 | 49821 | 49943 | 47510 |
| Hard Coal old | 27998 | 20183 | 6238 | 0 |
| Hard Coal new | 5226 | 3842 | 890 | 705 |
| CCGT new | 22627 | 16809 | 54881 | 48169 |

Looking forward, the RP-TM&CI model could be used to analyze the co-optimization of the water value in hydro storage with the storage value of short-term storage such as batteries. This kind of analysis could improve traditional hydrothermal dispatch models in which short-term storage is rarely considered. Moreover, the RP-TM&CI model could be extended to a stochastic model to consider uncertainty in renewable energy production or hydro inflows for long-term storage. Therefore, the main challenge in this topic is the representation at the same time of long- and short-term uncertainties, such as in [26].

APPENDIX

Tables III and IV summarize the main assumptions and results for the four different visions in the case study according to [22].

ACKNOWLEDGMENT

The authors would like to thank three anonymous referees for their comments and suggestions.

**Diego A. Tejada-Arango** (M'11) received the B.Sc. degree in electrical engineering from the "Universidad Nacional," Medellin, Colombia, in 2006, the M.Sc. degree in electrical engineering from the Universidad de Antioquia, Medellin, in 2013, with the GIMEL Group, and the Master's degree in research in engineering systems modeling from Universidad Pontificia Comillas, Madrid, Spain, in 2017. He is currently working toward the Ph.D. degree at the Instituto de Investigación Tecnológica, Universidad Pontificia Comillas de Madrid, Madrid. His research interests include transmission expansion planning, and planning and operation of electric energy systems.

**Maya Domeshek** received the B.A. degree in physics from Smith College, Northampton, MA, USA, in 2018. Her work on this project at the Universidad Pontificia Comillas de Madrid was funded by the WINDINSPIRE grant.

**Sonja Wogrin** (M'13) received the Dipl. Ing. (5-year) degree in technical mathematics from the Graz University of Technology, Graz, Austria, in 2008, the M.S. degree in computation for design and optimization from the Massachusetts Institute of Technology, Cambridge, MA, USA, in 2008, and the Ph.D. degree from the Instituto de Investigación Tecnológica, Universidad Pontificia Comillas de Madrid, Madrid, Spain, in 2013.

Since 2009, she has been a Researcher and an Assistant Professor, since 2013, with the Universidad Pontificia Comillas, Madrid. Her research interests include the area of decision support systems in the energy sector, optimization, and in particular the generation capacity expansion problem.

**Efraim Centeno** received the degree in industrial engineering and the Ph.D. degree in industrial engineering from the Universidad Pontificia Comillas de Madrid, Madrid, Spain, in 1991 and 1998, respectively. He belongs to the research staff with the Instituto de Investigación Tecnológica, Universidad Pontificia Comillas. His research interests include the planning and development of electric energy systems.